\documentclass[12pt]{article}

\usepackage{graphicx}
\usepackage{parskip}
\usepackage{siunitx}
\usepackage[utf8]{inputenc}
\usepackage{amsmath,amsfonts,amssymb,mathtools}
\usepackage{graphicx,float}
\usepackage{caption}
\usepackage{tikz}
\usetikzlibrary{automata, arrows.meta, positioning}

\usepackage[colorlinks,allcolors=black]{hyperref}
\newtheorem{theorem}{Theorem}

\newtheorem{proposition}[theorem]{Proposition}

\title{On tournaments' combinatorics}
\author{Joaquín Castañeda and Vinicio Gómez}

\date{}

\begin{document}

\maketitle

\begin{abstract}
\noindent \textit{In this paper we present a new proof of a proposition presented in the article: Complexes of tournaments, directionality filtrations and persistent homology [1]. This paper is part of Joaquín Castañeda's undergraduate thesis.
}\\


\end{abstract}

\noindent \textbf{Definition 1.1} [1] For a non-negative integer n, an n-tournament is a digraph with no reciprocal edges whose underlying undirected graph is an n-clique. An $n-$tournament is said to be

\begin{enumerate}
    \item \textit{transitive}, if its edge orientation defines a total ordering on its vertex set,
    \item \textit{regular}, if for each vertex $v \in \sigma$ its in-degree and out-degree are equal.
\end{enumerate}

\noindent \textbf{Definition 1.2} [1] Let $\mathcal{G} = (V,E)$ be a digraph. For a vertex $v \in V$, the following is defined:

\begin{enumerate}
    \item \textit{In-degree} of $v$ in $\mathcal{G}$, $i_{\mathcal{G}} (v)$, as the number of edges coming into $v$.
    \item \textit{Out-degree} of $v$ in $\mathcal{G}$, $o_{\mathcal{G}} (v)$, as the number of edges going out of $v$.
    \item \textit{Signed degree} of $v$ in $\mathcal{G}$ by  $sd_{\mathcal{G}}(v) = i_{\mathcal{G}}(v) - o_{\mathcal{G}}(v)$.
\end{enumerate}
    
\noindent Let $U \subseteq V$, a subset of vertices.

\begin{enumerate}
    \item The \textit{signed degree} of $U$ relative to $\mathcal{G}$ is defined as $sd_{\mathcal{G}}(U) = \sum_{v\in U} sd_{\mathcal{G}}(v)$.
    \item The \textit{directionality} of $U$ relative to$\mathcal{G}$ is defined as $Dr_\mathcal{G}(U) = \sum_{v\in U} sd_{\mathcal{G}}(v)^2$.
\end{enumerate}

\noindent \textbf{Definition 1.3} [1] For any \textit{$n$-tournament} $\sigma$ in a graph $\mathcal{G}$, let $V_{\sigma}$ denote the vertex set of $\sigma$ and define:

\begin{enumerate}
    \item \textit{Local directionality}: $Dr(\sigma) = Dr_{\sigma} (V_\sigma)$
    \item Let \textit{$c_3 (\sigma)$} denote the number of regular $3$-sub-tournaments in $\sigma$.
\end{enumerate}

\begin{proposition} Let $\sigma$ a tournament and $v,w \in \sigma$ two vertices such that the following requirement are met:

\begin{enumerate}
    \item $sd_{\sigma} (v) \geq sd_{\sigma}(w)$
    \item $i_{\sigma}(w) - i_{\sigma}(v) = k$ 
\end{enumerate}

Let $\tau$ the tournament formed by changing the direction of the edge between $v$ and $w$.

\textbf{First case:} $(v,w) \in \sigma$, then $Dr(\tau) = Dr(\sigma) - 8(k-1)$.

\textbf{Second case:} $(w,v) \in \sigma$, then $Dr(\tau) = Dr(\sigma) + 8(k+1)$.\\
\end{proposition}

\textbf{Proof:} For any vertex $u \in \sigma$ it is true that if $i_{\sigma}(u) + o_{\sigma}(u) = n - 1$, then $o_{\sigma}(u) = n - 1 - i_{\sigma}(u)$.\\

 \noindent Compute:
 
 \begin{center}
      $sd_{\sigma}(v) - sd_{\sigma}(w) = [i_{\sigma}(v) - o_{\sigma}(v)] - [i_{\sigma}(w) - o_{\sigma}(w)] = 2i_{\sigma}(v) - 2 i_{\sigma}(w)$\\ $= -2[i_{\sigma}(w) - i_{\sigma}(v)] = -2k$.
 \end{center}

 \noindent Let $sd_{\sigma}(v) = a$ and $sd_{\sigma}(w) = b$, thus $a - b = -2k$.\\
 
 \noindent \textbf{First case:} If the direction of the edge $(v,w)$ is changed, then for a vertex $u \in V_{\sigma} \setminus \{v,w\}$, $sd_{\sigma} (u) = sd{\tau} (u)$. Moreover,  $sd_{\tau}(v) = a + 2$ and $sd_{\tau}(w) = b-2$.\\
 
 \noindent Thus,
 
 \begin{center}
    $Dr(\tau) = Dr(\sigma) - a^2 - b^2 + (a+2)^2 + (b-2)^2 = Dr(\sigma) + 4(a-b) + 8 = Dr(\sigma) + 4(-2k) + 8 = Dr(\sigma) -8(k - 1)$.\\  
 \end{center}

\noindent \textbf{Second case:} As in the previous case, for a vertex $u \in V_{\sigma} \setminus \{v,w\}$, $sd_{\sigma} (u) = sd{\tau} (u)$ and for $v, w$ $sd_{\tau}(v) = a-2$ y  $sd_{\tau}(w) = b + 2$.\\

\noindent Thus, 

\begin{center}
     $Dr(\tau) = Dr(\sigma) - a^2 - b^2 + (a-2)^2 + (b+2)^2 = Dr(\sigma) + 4(b-a) + 8$\\  $= Dr(\sigma) + 8(k + 1)$.
\end{center}

\noindent \begin{proposition} Let $\sigma$ be a tournament that satisfies the same hypothesis as in the previous proposition.\\

\noindent Let $\tau$ the tournament formed by changing the direction of the edge between $v$ and $w$.\\

\textbf{First case:} $(v,w) \in \sigma$, then $ \tau$ has $(k-1)$ regular $3-$cycles more than $\sigma.$\\

\textbf{Second case:} $(w,v) \in \sigma$, then $ \tau$ has $(k+1)$  regular $3-$cycles less than $\sigma.$\\
\end{proposition}

\textbf{Proof:}
\textbf{First case:} Given $v, w$ and $u$ another vertex, then the $3-$cycle $\{u,v,w\} = \rho$ lies in one and only one of the following:

\begin{enumerate}
    \item It is regular in $\sigma$ but not in $\tau$.
    \item It is not regular in $\sigma$ and becomes regular in $\tau$.
    \item It is not regular in $\sigma$ and still not regular in $\tau$.
\end{enumerate}

\noindent If $\rho$ is such that $(u,v), (v,w) \in \rho$, then $\rho$ is type 1 or 2. The maximum number $D$ of regular $3-$cycles that could be destroyed is equal to the number of $3-$cycles of those two types, i.e the number of $3-$cycles in $\sigma$ with an edge incident on v. Such number $D$ is equal to $i_{\sigma}(v) = i_{\tau}(v)-1$.

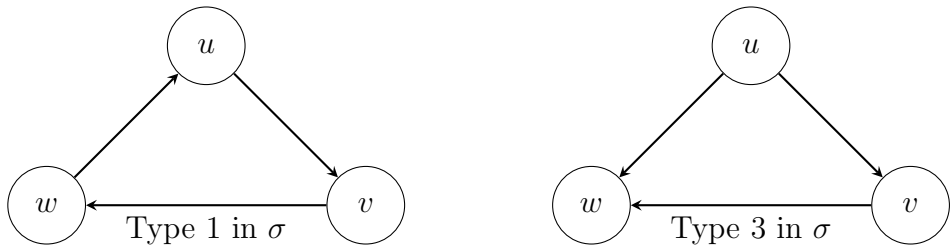
\begin{center}
\begin{tikzpicture} [node distance = 3cm, on grid, auto]
 
\node (v1) [state, initial text = {}] {$w$};
\node (v2) [state, above right = of v1] {$u$};
\node (v3) [state, below right = of v2] {$v$};
 
\path [-stealth, thick]
    (v1) edge node {}   (v2)
    (v2) edge node {}   (v3)
    (v3) edge node {Type 1 in $\sigma$}   (v1);

\node (v4) [state, right = of v3] {$w$};
\node (v5) [state, above right = of v4] {$u$};
\node (v6) [state, below right = of v5] {$v$};

\path[-stealth, thick]
    (v5) edge node {}   (v4)
    (v5) edge node {}   (v6)
    (v6) edge node {Type 3 in $\sigma$}   (v4);
\end{tikzpicture}
\end{center}

\noindent Now, the $3-$cycles $\rho$ in that could become regular are those formed by edges $(u,w)$ and $(v,w)$, but the edge between vertices $u$ and $v$ can go in any direction. So, the maximum number $C$ of $3-$cycles that could become regular is equal to the number of $3-$cycles of Types 2 and 3, i.e the number of $3-$cycles in $\tau$ with an edge incident on w. Such number $C$ is equal to $i_{\tau}(w) = i_{\sigma}(w)-1$.\\

\begin{center}
\begin{tikzpicture} [node distance = 3cm, on grid, auto]
 
\node (v1) [state, initial text = {}] {$w$};
\node (v2) [state, above right = of v1] {$u$};
\node (v3) [state, below right = of v2] {$v$};
 
\path [-stealth, thick]
    (v2) edge node {}   (v1)
    (v3) edge node {}   (v2)
    (v3) edge node {Type 2 in $\sigma$}   (v1);
    
\end{tikzpicture}
\end{center}

\noindent Thus, the number $C-D = i_{\sigma}(w) - i_{\sigma}(v) - 1 = k - 1$ represents the difference $c_3(\tau) - c_3(\sigma)$.\\

\noindent Therefore, there are $(k-1)$ regular $3-$cycles more in $\tau$ than in $\sigma$.\\

\noindent \textbf{Second case:} The proof for this case is similar to the previous case.\\

\noindent However, now $D = i_{\tau}(w) - 1 = i_{\sigma}(w)$ and $C = e_{\tau}(v) = i_{\sigma}(v) - 1 $. Then, the number $D-C = i_{\sigma}(w) - (i_{\sigma}(v) - 1) = k + 1$ represents the difference $c_3(\sigma) - c_3(\tau)$.\\

\noindent Therefore, there are $(k+1)$ regular $3-$cycles less in $\tau$ than in $\sigma$.

\begin{proposition} \text{[1]} Let $\sigma$ be an $n-$tournament. Then

\begin{center}
    $Dr(\sigma) = 2\binom{n+1}{3} - 8c_3(\sigma)$.
\end{center}

\end{proposition}

\textbf{Proof:}
\noindent Proceed by induction on $n$.\\

\noindent For $n=3$, $\sigma$ is regular or transitive. In the first case each vertex has signed degree equal to zero, then $Dr(\sigma) = 0$. If $\sigma$ is transitive, the signed degree of its vertices are $2, -2$ and $0$. Therefore $Dr(\sigma) = 8$. Thus the claim follows. Assuming the formula holds for all $n-$tournaments and prove it holds for $n$-tournaments.\\

\noindent \textbf{Particular case:} Let $\sigma$ an $n$-tournament. Add a new vertex $v$ and form a $(n+1)$-tournament $\tau$ where all the new edges are of the form $(v,u)$ or $(u,v)$.\\

\noindent \textbf{Observation:} All the $3$-cycles which contain the vertex $v$ are transitive. Thus, $c_3(\sigma) = c_3(\tau)$. Therefore, if the proposition is true

\begin{center}
  $Dr(\tau) - Dr(\sigma) =  2 \binom{n+2}{3} -  2 \binom{n+1}{3}$.
\end{center}

\noindent Note that $\binom{n+2}{3}$ is the number of $3-$cycles in an $(n+2)-$tournament $\rho$ and $\binom{n+1}{3}$ is the number of $3-$cycles in $\eta = \rho \setminus \{w\}$, where $w$ is any vertex. Then $\binom{n+2}{3} -  \binom{n+1}{3}$ is the number of $3-$cycles in $\rho \setminus \eta$. That is equal to the number of  edges in $\eta$. Thus $2\Big[\binom{n+2}{3} - \binom{n+1}{3} \Big]= 2\Big[\binom{n+1}{2}\Big] = n^2 + n$.\\

\noindent On the other hand, 

\begin{center}
    $Dr(\tau) - Dr(\sigma) = \sum_{u\in \tau} [i_{\tau}(u) - o_{\tau}(u)]^2 - \sum_{u\in \sigma} [i_{\sigma}(u) - o_{\sigma}(u)]^2.$
\end{center}

\noindent Regardless of the direction of the edges in $\tau \setminus \sigma$, $[i_{\tau}(v) - o_{\tau}(v)]^2 = n^2$. Then 

\begin{center}
    $Dr(\tau)= \sum_{u\in \sigma} [i_{\tau}(u) - o_{\tau}(u)]^2 + n^2$ and  $Dr(\tau) - Dr(\sigma) = \sum_{u\in \sigma} [i_{\tau}(u) - o_{\tau}(u)]^2 + n^2 - \sum_{u\in \sigma} [i_{\sigma}(u) - o_{\sigma}(u)]^2$\\
\end{center}

\noindent For the formula to be true, it is sufficient to prove that 

\begin{center}
    $\sum_{u\in \sigma} [i_{\tau}(u) - o_{\tau}(u)]^2 - \sum_{u\in \sigma} [i_{\sigma}(u) - o_{\sigma}(u)]^2 = \sum_{u\in \sigma} [i_{\tau}(u) - o_{\tau}(u)]^2 - [i_{\sigma}(u) - o_{\sigma}(u)]^2 = n$.
\end{center}

\noindent If the edges in $\tau \setminus \sigma$ come out of $v$ then, for each vertex $u \in \sigma$,  $o_{\tau}(u) = o_{\sigma}(u)$ and $i_{\tau}(u) = i_{\sigma}(u) + 1$. For the case where all the edges in $\tau \setminus \sigma$ incident on $v$ then, $o_{\tau}(u) = o_{\sigma}(u) + 1 $ y $i_{\tau}(u) = i_{\sigma}(u)$. In either case

\begin{center}
    $\sum_{u\in \sigma} [i_{\tau}(u) - o_{\tau}(u)]^2 - [i_{\sigma}(u) - o_{\sigma}(u)]^2 = \sum_{u\in \sigma} [1 + (i_{\sigma}(u) - o_{\sigma}(u))]^2 - [i_{\sigma}(u) - o_{\sigma}(u)]^2$\\
    
    $= \sum_{u\in \sigma} 1 + 2[e_{\sigma}(u) - s_{\sigma}(u)] = n + 2\Big( \sum_{u\in \sigma} e_{\sigma}(u) - \sum_{u\in \sigma} s_{\sigma}(u)\Big) = n + 2(0) = n$.
\end{center}

\noindent Therefore, the formula is valid for this particular case.\\

\noindent \textbf{General case:} Suppose that $\tau$ has $m$ edges that are incident on the vertex $v$. Consider the tournament $\rho$ where all edges go out from $v$ and changing, one by one, the direction of the appropriate edges to form $\tau$. With the help of the two previous propositions we have that:

\begin{center}
    $Dr(\tau) = Dr(\rho) - 8K$.\\
\end{center}

\noindent Where $K = c_3(\tau) - c_3(\rho)$. Therefore, 

\begin{center}
    $Dr(\tau) = 2 \binom{n+2}{3} - 8c_3(\rho) - 8[c_3(\tau) - c_3(\rho)] = Dr(\tau) = 2 \binom{n+2}{3} - 8c_3(\tau)$.
\end{center} 

\subsection*{Bibliography}
[1] D. Govc, R. Levi, J. P. Smith,Complexes of tournaments, directionality filtrations and persistent homology, 2021.

\end{document}